\documentclass[final]{siamart220329}   

\usepackage{pdflscape}
\usepackage{siampaper}
\usepackage[shortlabels]{enumitem}

\usepackage[
  createShortEnv,
  conf={one big link},
  commandRef=Cref]{proof-at-the-end}

\makeatletter
\def\thmt@innercounters{section,equation,theorem}
\makeatother

\NewDocumentEnvironment{proofEE}{O{}+b}{%
  \begin{proofE}[#1]
    #2
    \space
  \end{proofE}
}{}

\newcommand{\SYMMLQ}{\textsc{Symmlq}\xspace}
\newcommand{\CGM}{\textsc{Cg}\xspace}
\newcommand{\CR}{\textsc{Cr}\xspace}
\newcommand{\CAR}{\textsc{cAr}\xspace}
\newcommand{\MINRES}{\textsc{Minres}\xspace}
\newcommand{\MINARES}{\textsc{MinAres}\xspace}

\newcommand{\MINRESQLP}{\textsc{Minres-qlp}\xspace}

\newcommand{\LSQR}{\textsc{Lsqr}\xspace}
\newcommand{\LSMR}{\textsc{Lsmr}\xspace}

\def\T{^T\!}
\def\TT{^{\!\!\!T}\!}

\newlength{\forwidth}
\makeatletter
\g@addto@macro\bfseries{\boldmath}
\makeatother

\newcommand{\bmat}[1]{\begin{bmatrix} #1 \end{bmatrix}}

\newcommand{\smallbmat}[1]{\hbox{\scriptsize $\bmat{#1}$}}

\newcommand{\norm}[1]{\|#1\|}
\newcommand{\Norm}[1]{\left\|#1\right\|}

\newcommand{\bb}[1]{\bar{\bar{#1}}}
\newcommand{\mr}[1]{\mathring{#1}}

\newcommand{\gammabar}{\bar\gamma}
\newcommand{\lambdabar}{\bar\lambda}
\newcommand{\mubar}{\skew3\bar \mu}
\newcommand{\phibar}{\bar{\phi}}

\newcommand{\zbar}{\skew{2.8}\bar z}
\newcommand{\zetabar}{\bar{\zeta}}
\newcommand{\psibar}{\bar{\psi}}
\newcommand{\thetabar}{\bar{\theta}}
\newcommand{\taubar}{\bar{\tau}}
\newcommand{\chibar}{\bar{\chi}}
\newcommand{\pibar}{\bar{\pi}}

\newcommand{\chat}{\hat{c}}
\newcommand{\shat}{\hat{s}}

\newcommand{\gammahat}{\skew4\hat\gamma}
\newcommand{\lambdahat}{\hat\lambda}
\newcommand{\Lhat}{\hat{L}}
\newcommand{\Phat}{\hat{P}}

\newcommand{\ctilde}{\tilde c}
\newcommand{\Qtilde}{\widetilde Q}
\newcommand{\stilde}{\tilde s}

\newcommand{\xstar}{x^{\star}}

\usepackage[utf8]{inputenc}
\usepackage[T1]{fontenc}
\usepackage{booktabs}

\usepackage{newfloat}
\DeclareFloatingEnvironment[
  name=TripleAlgorithm,
  placement=ht,
  fileext=lod]{triplealgorithm}

\graphicspath{{graphics/}}  
\usepackage{tikzscale}
\usepackage{pgfplots}
\pgfplotsset{compat=newest}
\usetikzlibrary{backgrounds}
\usetikzlibrary{intersections}
\usetikzlibrary{external}
\usetikzlibrary{math}

\iffalse
  \newcommand*{\includetikzgraphics}[2][]{%
    \tikzsetnextfilename{#2}%
    \includegraphics[#1]{#2.tikz}
  }
\else
  \newcommand*{\includetikzgraphics}[2][]{%
    \includegraphics[#1]{#2.pdf}
  }
\fi

\makeatletter
\renewcommand{\todo}[2][]{\tikzexternaldisable\@todo[#1]{#2}\tikzexternalenable}
\makeatother

\newcommand{\cahiernumber}{40}  

\usepackage{kbordermatrix}

\title{%
  {\MINARES}:\@ An Iterative Solver for \\ Symmetric Linear Systems
}

\author{%
  Alexis Montoison%
  \thanks{%
    GERAD and Department of Mathematics and Industrial Engineering,
    Polytechnique Montr\'eal, QC, Canada.
    E-mail: \mailto{alexis.montoison@polymtl.ca}.
    Research supported by an FRQNT grant and an excellence scholarship of the IVADO institute.
  }
  \and
  Dominique Orban%
  \thanks{%
    GERAD and Department of Mathematics and Industrial Engineering,
    Polytechnique Montr\'eal, QC, Canada.
    E-mail: \mailto{dominique.orban@gerad.ca}.
    Research partially supported by an NSERC Discovery Grant.
  }
  \and
  Michael A. Saunders%
  \thanks{%
    Systems Optimization Laboratory,
    Department of Management Science and Engineering, Stanford University, Stanford, CA, USA.
    E-mail: \mailto{saunders@stanford.edu}.
    \hfill Version of \today.%
  }
}
\date{\today}
\headers{\MINARES}{A. Montoison, D. Orban, and M. A. Saunders}

\begin{document}

  \maketitle

  \thispagestyle{firstpage}
  \pagestyle{myheadings}

  \begin{abstract}
    We introduce an iterative solver named \MINARES for symmetric linear systems $Ax \approx b$, where $A$ is possibly singular.
    \MINARES is based on the symmetric Lanczos process, like \MINRES and \MINRESQLP, but it minimizes $\norm{Ar_k}$ in each Krylov subspace rather than $\norm{r_k}$, where $r_k$ is the current residual vector.
    When $A$ is symmetric, \MINARES minimizes the same quantity $\norm{Ar_k}$ as \LSMR, but in more relevant Krylov subspaces, and it requires only one matrix-vector product $Av$ per iteration, whereas \LSMR would need two.
    Our numerical experiments with \MINRESQLP and \LSMR show that \MINARES is a pertinent alternative on consistent symmetric systems and the most suitable Krylov method for inconsistent symmetric systems.
    We derive properties of \MINARES from an equivalent solver named \CAR that is to \MINARES as \CR is to \MINRES, is not based on the Lanczos process, and minimizes $\norm{Ar_k}$ in the same Krylov subspace as \MINARES. 
    We establish that \MINARES and \CAR generate monotonic $\|x_k - \xstar\|$, $\|x_k - \xstar\|_A$ and $\|r_k\|$ when $A$ is positive definite.
  \end{abstract}




\section{Introduction}

Suppose $A \in \R^{n \times n}$ is a large symmetric matrix for which matrix-vector products $Av$ can be computed efficiently for any vector $v \in \R^n$.
We present a Krylov subspace method called \MINARES for computing a solution to the following problems:
\begin{alignat}{2}
\label{eq:Ax=b}
  & \text{Symmetric linear systems:} && \quad Ax = b,
\\
\label{eq:ls}
& \text{Symmetric least-squares problems:} && \quad \min \norm{Ax - b},
\\
\label{eq:nullvector}
& \text{Symmetric nullspace problems:}  && \quad Ar = 0,
\\
\label{eq:eigenvector}
& \text{Symmetric eigenvalue problems:} && \quad Ar = \lambda r,
\\[-4pt]
\label{eq:singularvector}
& \text{Singular value problems for rectangular $B$:} && \quad      \bmat{& B \\ B\T & } \bmat{u \\ v} = \sigma \bmat{u \\ v}.
\end{alignat}
If $A$ is nonsingular, problems~\eqref{eq:Ax=b}--\eqref{eq:ls}
have a unique solution $\xstar$.  When $A$ is singular, if $b$ is not in the
range of $A$ then \eqref{eq:Ax=b} has no solution; otherwise, \eqref{eq:Ax=b}--\eqref{eq:ls} have an infinite
number of solutions, and we seek the unique solution $\xstar$
that minimizes $\norm{x}$.
Whenever $\xstar$ exists, it solves the problem
\begin{equation}
    \label{eq:AAx=Ab}
    \min \tfrac{1}{2} \|x\|^2 \quad \st \quad A^2 x = Ab.
\end{equation}

Let $x_k$ be an approximation to $\xstar$
with residual $r_k = b - Ax_k$.
If $A$ were unsymmetric or rectangular, applicable solvers for~\eqref{eq:Ax=b}--\eqref{eq:ls} would be \LSQR \cite{paige-saunders-1982} and \LSMR \cite{fong-saunders-2011}, which reduce
$\norm{r_k}$ and $\norm{A\T r_k}$ respectively within the $k$th Krylov subspace $\mathcal{K}_k(A\T A, A\T b)$ generated by the Golub-Kahan bidiagonalization on $(A,b)$ \cite{golub-kahan-1965}.

For~\eqref{eq:Ax=b}--\eqref{eq:singularvector}, we propose an algorithm \MINARES that solves~\eqref{eq:AAx=Ab} by reducing $\norm{Ar_k}$ within the $k$th Krylov subspace $\mathcal{K}_k(A,b)$ generated by the symmetric Lanczos process on $(A,b)$ \cite{lanczos-1950}.
Thus when $A$ is symmetric, \MINARES minimizes the same quantity $\norm{Ar_k}$ as \LSMR, but in different (more effective) subspaces, and it requires only one matrix-vector product $Av$ per iteration, whereas \LSMR would need two.

Qualitatively, certain residual norms decrease smoothly for these iterative methods, but other norms are more erratic as they approach zero.  It is ideal if stopping criteria involve the smooth quantities. 
For \LSQR and \LSMR on general (possibly rectangular) systems, $\norm{r_k}$ decreases smoothly for both methods.
We observe that while \LSQR is always ahead by construction, it is never by very much.  Thus on consistent systems $Ax=b$, \LSQR 
may terminate slightly sooner.
On inconsistent systems $Ax \approx b$,
the comparison is more striking.
$\norm{A\T r_k}$ decreases erratically for \LSQR but smoothly for \LSMR,
and there is usually a significance difference between the two.
Thus \LSMR may terminate significantly sooner \citep{fong-saunders-2011}.

Similarly for \MINRES \citep{paige-saunders-1975} and \MINARES, $\norm{r_k}$ decreases smoothly for both methods,
and on consistent symmetric systems $Ax=b$, \MINRES may have a small advantage.
On inconsistent symmetric systems $Ax \approx b$,
$\norm{Ar_k}$ decreases erratically for \MINRES and its variant \MINRESQLP \citep{choi-paige-minres-2011} but smoothly for \MINARES,
and there is usually a significant difference between them.
Thus \MINARES may terminate sooner.

We introduce \CAR, a new conjugate direction method similar to \CGM and \CR and equivalent to \MINARES when $A$ is SPD.
We prove that $\norm{r_k}$, $\norm{x_k - \xstar}$ and $\norm{x_k - \xstar}_A$ decrease monotonically for \CAR and hence \MINARES when $A$ is positive definite.

\subsection{Notation}
A symmetric positive definite matrix is said to be SPD.
For a vector $v_k$, $\norm{v_k}$ denotes the Euclidean norm of $v_k$, and for an SPD matrix \(A\), the
$A$-norm of $v_k$ is $\norm{v_k}_A^2 = v\T A v$.
For a matrix $V_k$, $\norm{V_k}$ may be any norm.
Vector $e_j$ is the $j$th column of an identity matrix $I_k$ of size dictated by the context.
An approximate solution $x_k$ has residual $r_k = b - Ax_k$, and
$\xstar$ is the unique solution of $Ax = b$ if $A$ is nonsingular, or the minimum-norm solution of $A^2x = Ab$ otherwise.
$\mathcal{K}_k(A, b)$ is the \emph{Krylov subspace} $\{b, Ab, \dots, A^{k-1}b\}$.
We abusively write $z = (\zeta_1, \dots, \zeta_n)$ to represent the column vector $z = \bmat{\zeta_1 & \dots & \zeta_n}\T$.
If $H$ is SPD and $\{d_1, \dots, d_k\}$ is a set of non-zero vectors, the vectors are $H$-conjugate if $d_i\T H d_j = 0$ for $i \ne j$.
If $H = I$, conjugacy is equivalent to the usual notion of orthogonality.

\section{Applications}
 
\subsection{Null vector, eigenvector, and singular value problems}

Given a symmetric $A$ and nonzero $b$, \MINARES solves $A^2 x = Ab$ even if $A$ is singular. If $b$ is random and $A$ is singular, $r = b - Ax$ is unlikely to be zero, but it will be a nonzero nullvector of $A$ because $Ar = 0$.

If an eigenvalue $\lambda$ of $A$ is known, we can use it as a shift in the Lanczos process with a random starting vector $b$ to find a null vector $r$ such that $(A - \lambda I)r = 0$. Then $r$ is an eigenvector because $Ar = \lambda r$.
\MINARES is effectively implementing the inverse power method \citep{trefethen-bau-1997,golub-vanloan-2013} to obtain the eigenvector in one iteration.  If $\lambda$ is approximate, \MINARES can implement Rayleigh quotient iteration \citep{trefethen-bau-1997,golub-vanloan-2013} to obtain increasingly accurate eigenpair estimates.

Similarly, if a singular value $\sigma$ is known for a rectangular matrix $B$, the singular value problem $Bv = \sigma u$, $B\T u = \sigma v$ may be reformulated as a null vector problem or eigenvalue problem:
$$\left(\bmat{& B \\ B^T &} - \sigma I\right) \bmat{u \\ v} = 0
  \quad \Longleftrightarrow \quad
  \bmat{ & B \\ B^T & } \bmat{u \\ v} = \sigma \bmat{u \\ v},
$$
for which \MINARES may be used to implement inverse iteration or Rayleigh quotient iteration
(although an algorithm based on the Golub-Kahan bidiagonalization of $B$ would be preferable).



\subsection{Singular systems with semi-positive definite matrices}

Inconsistent (singular) symmetric systems could arise from discretized semidefinite Neumann boundary value problems \citep[sect.~4]{kaasschieter-1988}.
Measurement errors will be random, so $b$ is unlikely to be in the range of singular $A$.

Another potential application is large, singular, symmetric, indefinite Toeplitz least-squares problems as described in \citep[sec.~5]{gallivan-thirumalai-dooren-vermaut-1996}.
Rank-deficient Toeplitz matrices arise in image reconstruction and system identification problems.
In both cases, $A$ is a semi-positive definite matrix and \MINARES is a suitable solver.

\section{Symmetric systems} \label{sec:symmetricA}

\begin{algorithm}[t]
  \caption{Lanczos process}
  \label{alg:lanczos}
  \begin{algorithmic}[1]
    \Require $A$, $b$
    \State $v_0 = 0$
    \State $\beta_1 v_1 = b$  \Comment{$\beta_1 > 0$ so that $\norm{v_1} = 1$}
    \For{$k=1,2,\dots$}
      \State $q_k = Av_k - \beta_k v_{k-1}$
      \State $\alpha_k = v_k^T q_k$
      \State $q_k = q_k - \alpha_k v_k$
      \State $\beta_{k+1} = \|q_k\|$
      \If{$\beta_{k+1} = 0$}
          \State $\ell = k$; \ \textbf{return} $\ell$
      \Else
          \State $v_{k+1} = q_k / \beta_{k+1}$ \Comment{$\beta_{k+1} > 0$ so that $\norm{v_{k+1}} = 1$}
      \EndIf
    \EndFor
  \end{algorithmic}
\end{algorithm}

With $A$ symmetric and starting vector $b$, we make use of the symmetric Lanczos process \cite{lanczos-1950} of \cref{alg:lanczos}.
After \(k\) iterations the situation may be summarized as
\begin{subequations}
    \begin{alignat}{2}
    A V_k &= V_k T_k + \beta_{k+1} v_{k+1} e_k^T &&= V_{k+1} T_{k+1,k},  \label{eq:lanczos1}
  \\ V_k^T V_k &= I_k,  \label{eq:lanczos2}
  \end{alignat}
\end{subequations}
where
\vspace*{-10pt}
\begin{equation*}
  V_k := \bmat{v_1 & \dots & v_k},
  \qquad T_k =
  \smallbmat{\alpha_1 & \beta_2
     \\ \beta_2  & \alpha_2 & \ddots
     \\          & \ddots   & \ddots   & \beta_k
     \\          &          & \beta_k  & \alpha_k},
  \qquad T_{k+1,k} = \bmat{T_{k} \\ \beta_{k+1} e_{k}^T}.
\end{equation*}
In exact arithmetic, $V_k$ is an orthonormal basis of $\mathcal{K}_k(A, b)$.  The Lanczos process terminates after $\ell \le n$ iterations when $\beta_{\ell+1}=0$, and we then have $A V_{\ell} = V_{\ell} T_{\ell}$, where square $T_{\ell}$ is nonsingular if and only if $b \in \text{range}(A)$ \citep[sec.~2.1 property 4]{choi-paige-minres-2011}.
$T_{k+1,k}$ has full column rank $k$ for all $k < \ell$ \citep[sec.~2.1 property 2]{choi-paige-minres-2011} and the rank of $T_{\ell}$ is $\ell$ or $\ell-1$ but no less (because the first $\ell-1$ columns of $T_{\ell}$ are independent).

In finite arithmetic, \eqref{eq:lanczos1} holds to machine precision.
Reorthogonalization would be needed for \eqref{eq:lanczos2} to hold accurately,
but it is enough to note that we always have $\norm{V_k} = O(1)$.

\subsection{\CGM, \SYMMLQ, \MINRES, \MINARES}

As with \CGM \citep{hestenes-stiefel-1952},
 \SYMMLQ \citep{paige-saunders-1975}, and
 \MINRES \citep{paige-saunders-1975}, the goal of \MINARES is to solve symmetric problems $Ax \approx b$.
All methods define an approximate solution $x_k = V_k y_k$ at iteration $k$ (where $y_k$ is different for each method).
\MINARES chooses $y_k$ to minimize $\norm{Ar_k}$ in $\mathcal{K}_k(A, b)$, so that $\norm{Ar_k}$ is monotonically decreasing towards zero.
\MINARES is therefore well suited to singular inconsistent symmetric systems.
This case is difficult for the other methods because 
$\norm{x_k - \xstar}_A$, $\norm{x_k - \xstar}$ and $\norm{r_k}$ do not converge to zero and they are the quantities minimized respectively by \CGM, \SYMMLQ, and both \MINRES and \MINRESQLP.

\section{Derivation of \MINARES}

\subsection{Subproblems of \MINARES}

From \Cref{alg:lanczos} we have $Ab = \beta_1 \alpha_1 v_1 + \beta_1 \beta_2 v_2$ because $\beta_2 v_2 = Av_1 - \alpha_1 v_1$.
Hence
\begin{subequations} \label{eq:Ar_variants}
\begin{align}
  A r_k &= A (b - A V_k y_k)           \nonumber
\\      &= Ab - AV_{k+1} T_{k+1,k} y_k \nonumber
\\      &= \beta_1 \alpha_1 v_1 + \beta_1 \beta_2 v_2
           - V_{k+2} T_{k+2,k+1} T_{k+1,k} y_k \nonumber
\\      &= V_{k+2}(\beta_1 \alpha_1 e_1 + \beta_1 \beta_2 e_2
           - T_{k+2,k+1} T_{k+1,k} y_k),
           \quad k \le \ell - 2,  \label{eq:Ar}
\\ A r_{\ell-1} &= V_{\ell} (\beta_1 \alpha_1 e_1 + \beta_1   \beta_2 e_2 -  T_{\ell} T_{\ell,\ell-1} y_{\ell-1}), \label{eq:Ar2}
\\ A r_{\ell} &= V_{\ell} (\beta_1 \alpha_1 e_1 + \beta_1 \beta_2 e_2 - T_{\ell}^2 y_{\ell}). \label{eq:Ar3}
\end{align}
\end{subequations}
Theoretically, $V_{k}$ has orthonormal columns ($1 \le k \le \ell$), so that $\norm{x_k} = \norm{y_k}$ and 
$\norm{Ar_k}$ is minimized with $\|x_k\|$ of minimal norm if we define $y_k$ as the unique solution of the following subproblems:
\begin{subequations} \label{eq:sub-minares}
  \begin{align}
    \label{eq:sub1}
    \minimize{y_k \in \R^k} \quad & \norm{T_{k+2,k+1} T_{k+1,k} y_k - \beta_1 \alpha_1 e_1 - \beta_1 \beta_2 e_2}, \quad k \le \ell-2,
  \\
    \label{eq:sub2}
    \minimize{y_{\ell-1} \in \R^{\ell-1}} \quad & \norm{T_{\ell} T_{\ell,\ell-1} y_{\ell-1} - \beta_1 \alpha_1 e_1 - \beta_1 \beta_2 e_2},
  \\
    \label{eq:sub3}
    \minimize{y_{\ell} \in \R^{\ell}} \quad & \|y_{\ell}\|^2 \quad \st \quad T_{\ell}^2 y_{\ell} = \beta_1 \alpha_1 e_1 + \beta_1 \beta_2 e_2.
  \end{align}
\end{subequations}
We define $y_k$ from these subproblems even though $V_{k}$ does not remain orthonormal numerically.
In practice, we expect $\norm{Ar_k} \le \norm{Ar_{k-1}}$ unless $k$ becomes too large.

To be sure that the subproblems have unique solutions, we need to verify that
$T_{k+2,k+1} T_{k+1,k}$ has rank $k$ ($k \le \ell-2$),
$T_{\ell} T_{\ell,\ell-1}$ has rank $\ell-1$,
and $T_{\ell}^2 y_{\ell} = \beta_1 \alpha_1 e_1 + \beta_1 \beta_2 e_2$ is consistent even if $T_{\ell}$ is singular.
These results are proved in \Cref{theorem:TT_full_rank1}, \Cref{theorem:TT_full_rank2} and \Cref{theorem:Ar=0}.

\begin{theoremE}[][end, restate]
  \label{theorem:TT_full_rank1}
  For $k \le \ell-2$, $T_{k+2,k+1} T_{k+1,k}$ has rank $k$.
\end{theoremE}

\begin{proofEE}
From~\eqref{eq:Nk1} and~\eqref{eq:Nk} we have
$$
T_{k+2,k+1} T_{k+1,k} =\! \bmat{R_k\T R_k \\ (\varepsilon_{k-1} e_{k-1}\T + \gamma_k e_k\T) R_k \\ \varepsilon_k e_k\T R_k},
$$
where $R_k\T R_k$ has rank $k$ because $T_{k+1,k}$ and hence $R_k$ have full column rank.
\end{proofEE}

\begin{theoremE}[][end, restate]
  \label{theorem:TT_full_rank2}
  $T_{\ell} T_{\ell,\ell-1}$ has rank $\ell-1$.
\end{theoremE}

\begin{proofEE}
From~\eqref{eq:Nk2} and~\eqref{eq:Nk} we have
$$
T_{\ell} T_{\ell,\ell-1} =\! \bmat{R_{\ell-1}\T R_{\ell-1} \\ (\varepsilon_{\ell-1} e_{\ell-1}\T + \gamma_{\ell} e_{\ell}\T) R_{\ell-1}}\!,
$$
where $R_{\ell-1}\T R_{\ell-1}$ has rank $\ell-1$ because $T_{\ell,\ell-1}$ and $R_{\ell-1}$ have full column rank.
\end{proofEE}

\begin{theoremE}[][end, restate]
  \label{theorem:Ar=0}
  $T_{\ell}^2 y_{\ell} = \beta_1 \alpha_1 e_1 + \beta_1 \beta_2 e_2$ is consistent even if $T_{\ell}$ is singular.
\end{theoremE}

\begin{proofEE}
If $T_{\ell}$ is singular, the symmetry of $T_{\ell}$ and its complete orthogonal decomposition give
\begin{equation*}
T_{\ell} = Q \bmat{L & 0 \\ 0 & 0} P = P\T \bmat{L\T & 0 \\ 0 & 0} Q\T \quad \text{and} \quad T_{\ell}^2 = P\T \bmat{L\T L & 0 \\ 0 & 0} P,
\end{equation*}
where $Q$ and $P$ are orthogonal and rank$(L) = \ell-1$.
Thus,
\begin{align*}
   T_{\ell}^2 y_{\ell} - \beta_1 \alpha_1 e_1 - \beta_1 \beta_2 e_2 &=
   T_{\ell}^2 y_{\ell} - \beta_1 T_{\ell} e_1\\
   &=
      P\T
      \left(
         \bmat{L\T L & 0 \\ 0 & 0} P y_{\ell}
       - \beta_1 \bmat{L^T & 0 \\ 0 & 0} Q\T e_1
      \right)
\\ &= P\T \bmat{L\T L t_{\ell-1} - L\T u_{\ell-1} \\ 0},
\end{align*}
where $t_{\ell-1}$ and $u_{\ell-1}$ are the first $\ell-1$ components of $Py_{\ell}$ and $\beta_1 Q\T e_1$.
Because $L$ has full rank, $L\T L t_{\ell-1} = L\T u_{\ell-1}$ has a unique solution.
Then, $y_\ell = P^T \smallbmat{t_{\ell-1} \\ \omega}$ is a solution of $T_{\ell}^2 y_{\ell} = \beta_1 \alpha_1 e_1 + \beta_1 \beta_2 e_2$ for any $\omega$, which means the system is consistent.
\end{proofEE}

From~\eqref{eq:Ar3} and \Cref{theorem:Ar=0}, $Ar_{\ell} = V_{\ell}(T_{\ell}^2 y_{\ell} - \beta_1 \alpha_1 e_1 - \beta_1 \beta_2 e_2) = 0$.
Hence with definition \eqref{eq:sub3} we can conclude that $x_{\ell}$ is the solution $\xstar$ of \eqref{eq:AAx=Ab}.

\subsection{QR factorization of \texorpdfstring{$T_k$}{Tk}}

To solve~\eqref{eq:sub-minares}, we first need the QR factorization used by \MINRES:
\vspace*{-10pt}
\begin{equation}
  \label{eq:qr_Tk}
  T_{k+1,k} = Q_k \bmat{R_k \\ 0},
  \quad
  R_k =
  \smallbmat{
  \lambda_1  & \gamma_1  & \varepsilon_1 &        &                   \\
             & \lambda_2 & \gamma_2      & \ddots &                   \\
             &           & \lambda_3     & \ddots & \varepsilon_{k-2} \\
             &           &               & \ddots & \gamma_{k-1}      \\
             &           &               &        & \lambda_{k}
  },
\end{equation}
where $Q_k^T = Q_{k+1,k}\dots Q_{3,2}Q_{2,1}$ is an orthogonal matrix defined as a product of $2 \times 2$ reflections with the structure
\begin{equation*}
  \label{eq:reflection_structure}
  Q_{i+1,i} =
  \hbox{\scriptsize $
  \kbordermatrix{
         & 1 & \hdots & i-1 & i   & i+1            & i+2 & \hdots & k \\
  1      & 1 &        &     &     &                &     &        &   \\
  \vdots &   & \ddots &     &     &                &     &        &   \\
  i-1    &   &        & 1   &     &                &     &        &   \\
  i      &   &        &     & c_i & \phantom{-}s_i &     &        &   \\
  i+1    &   &        &     & s_i &          - c_i &     &        &   \\
  i+2    &   &        &     &     &                & 1   &        &   \\
  \vdots &   &        &     &     &                &     & \ddots &   \\
  k      &   &        &     &     &                &     &      & 1   \\
  }$}.
\end{equation*}
If we initialize \(Q_0 := I\), \(\lambdabar_1 := \alpha_1\), \(\gammabar_1 := \beta_2\), individual factorization steps may be represented as an application of $Q_{k+1,k}$ to $Q_{k-1}^T T_{k+1,k}$:
\[
  \kbordermatrix{
      & k   & k+1             \\
  k   & c_k & \phantom{-} s_k \\
  k+1 & s_k &          -  c_k
  }
  \!
  \kbordermatrix{
  & k              &        & k+1          & k+2         \\
  & \lambdabar_{k} & \vrule & \gammabar_k  & 0           \\
  & \beta_{k+1}    & \vrule & \alpha_{k+1} & \beta_{k+2}
  }
  =
  \!\!\!\!\!
  \kbordermatrix{
  & k                &        & k+1              & k+2             \\
  & \lambda_{k}      & \vrule & \gamma_k         & \varepsilon_{k} \\
  & 0                & \vrule & \lambdabar_{k+1} & \gammabar_{k+1}
  }.
\]

The reflection $Q_{k+1,k}$ zeroes $\beta_{k+1}$ on the subdiagonal of $T_{k+1,k}$ and affects three columns and two rows.
It is defined by
\begin{equation}
  \lambda_{k} = \sqrt{\lambdabar_{k}^2 + \beta_{k+1}^2},
  \quad
  c_k = \lambdabar_{k} / \lambda_{k},
  \quad
  s_k = \beta_{k+1} / \lambda_{k},
\end{equation}
and yields the following recursion for $k \ge 1$:
\begin{subequations}
  \begin{align}
  \label{qr1_formula}
   \gamma_{k}       &= c_k \gammabar_{k} + s_k
             \alpha_{k+1}, 
\\ \lambdabar_{k+1} &= s_k \gammabar_{k} - c_k
             \alpha_{k+1}, 
\\ \varepsilon_{k}  &= \phantom{-}s_{k} \beta_{k+2}, 
\\ \gammabar_{k+1}  &= - c_{k} \beta_{k+2}.
  \end{align}
\end{subequations}

\subsection{Definition of \texorpdfstring{$N_k$}{Nk}}

Let us define
\begin{subequations}
  \label{eq:Nk_definition}
  \begin{alignat}{3}
  N_{k} &:= T_{k+2,k+1} Q_k \bmat{I_k \\ 0}, \quad
        & \text{where}
        & \quad N_{k} R_k =
                   T_{k+2,k+1} T_{k+1,k}, \quad k \le \ell-2,
  \label{eq:Nk1} \\
  N_{\ell-1} &:= T_{\ell,\ell-1} Q_{\ell-1}
                    \bmat{I_{\ell-1} \\ 0}, \quad
           & \text{where}
           & \quad N_{\ell-1} R_{\ell-1} = T_{\ell} T_{\ell,\ell-1},
  \label{eq:Nk2} \\
  N_{\ell} &:= T_{\ell} Q_{\ell}, \quad
           & \text{where}
           & \quad N_{\ell} R_{\ell} = T_{\ell}^2.
  \label{eq:Nk3}
  \end{alignat}
\end{subequations}
Because $Q_k = Q_{2,1} Q_{3,2} \dots Q_{k+1,k}$, we have
\begin{subequations}
  \label{eq:ekQk}
  \begin{align}
  e_{k}^T Q_k & =
  e_{k}^T Q_{k,k-1} Q_{k+1,k}=
  s_{k-1} e_{k-1}^T - c_{k-1} c_{k} e_{k}^T - c_{k-1} s_k e_{k+1}^T,
  \\
  e_{k+1}^T Q_k & =
  e_{k+1}^T Q_{k+1,k} =
  s_{k} e_{k}^T - c_{k} e_{k+1}^T.
  \end{align}
\end{subequations}
Moreover, $T_{k+2,k+1} = \smallbmat{T_{k+1,k}^T \\ \beta_{k+1} e_k^T + \alpha_{k+1} e_{k+1}^T \\ \beta_{k+2} e_{k+1}^T}$ and the product $T_{k+2,k+1} Q_k$ can be determined in three parts.
From~\eqref{eq:qr_Tk}, $T_{k+1,k}^T Q_k = \left(Q_k^T T_{k+1,k}\right)^T = \bmat{R_k^T & 0}$,
and from~\eqref{eq:ekQk} we have
\begin{equation*}
  \begin{aligned}
    (\beta_{k+1} e_k^T + \alpha_{k+1} e_{k+1}^T) Q_k &= \beta_{k+1} s_{k-1} e_{k-1}^T
    + (\alpha_{k+1} s_k - \beta_{k+1} c_{k-1} s_k)e_k^T
 \\ &\hspace*{74pt}
    - (\alpha_{k+1}c_k + \beta_{k+1} c_{k-1} s_k) e_{k+1}^T
 \\ &= \varepsilon_{k-1}e_{k-1}^T + \gamma_k e_k^T - (\alpha_{k+1}c_k + \beta_{k+1} c_{k-1} s_k) e_{k+1}^T,
 \\ \beta_{k+2} e_{k+1}^T Q_k &= s_{k} \beta_{k+2} e_{k}^T - c_{k} \beta_{k+2} e_{k+1}^T
 \\ &= \varepsilon_k e_k^T - c_{k} \beta_{k+2} e_{k+1}^T.
  \end{aligned}
\end{equation*}
Thus, for $k \le \ell-2$ we obtain
\begin{equation}
  \label{eq:Nk}
  N_{k} = \bmat{R_k^T
                \\ \varepsilon_{k-1}e_{k-1}^T + \gamma_k e_k^T
                \\ \varepsilon_{k}e_{k}^T},
          \quad N_{\ell-1} = \bmat{R_{\ell-1}^T
                              \\ \varepsilon_{\ell-1}e_{\ell-1}^T + \gamma_{\ell} e_{\ell}^T},
          \quad
        N_{\ell} = R_{\ell}^T.
\end{equation}

\subsection{QR factorization of \texorpdfstring{$N_k$}{Nk}}

\vspace*{-10pt}
\begin{equation}
  \label{eq:qr_Nk}
  N_k = \Qtilde_k \bmat{U_k \\ 0},
  \quad
  U_k = \bmat{ 
    \mu_1 & \phi_1 & \rho_1 &        &
 \\       & \mu_2  & \phi_2 & \ddots &
 \\       &        & \mu_3  & \ddots & \rho_{k-2}
 \\       &        &        & \ddots & \phi_{k-1}
 \\ \phantom{\ddots} & &    &        & \mu_{k}
  },
\end{equation}
where $\Qtilde_k^T = \Qtilde_{k+2,k} \Qtilde_{k+1,k} \dots \Qtilde_{3,1}\Qtilde_{2,1}$ for $k\le \ell-2$ , and $\Qtilde_{\ell}^T = \Qtilde_{\ell-1}^T = \Qtilde_{\ell,\ell-1}\Qtilde_{\ell-2}^T$ are orthogonal matrices defined as a product of reflections.
If we initialize \(\mubar_1 := \lambda_1\), \(\gammahat_{1} := \gamma_1\) and \(\lambdahat_2 := \lambda_2\), individual factorization steps may be represented as an application of $\Qtilde_{k+1,k}$ to $\Qtilde_{k-1}^T N_{k}$:
\[
\kbordermatrix{
& k                  & k+1                        & k+2 \\
k   & \ctilde_{2k-1} & \phantom{-} \stilde_{2k-1} &     \\
k+1 & \stilde_{2k-1} &          -  \ctilde_{2k-1} &     \\
k+2 &                &                            & 1
}
\!
\kbordermatrix{
& k               &        & k+1              & k+2           \\
& \mubar_{k}      & \vrule &                  &               \\
& \gammahat_{k}   & \vrule & \lambdahat_{k+1} &               \\
& \varepsilon_{k} & \vrule & \gamma_{k+1}     & \lambda_{k+2}
}
=
\!\!\!\!\!
\kbordermatrix{
& k               &        & k+1          & k+2           \\
& \bb{\mu}_{k}    & \vrule & \phibar_{k}  &               \\
&                 & \vrule & \mubar_{k+1} &               \\
& \varepsilon_{k} & \vrule & \gamma_{k+1} & \lambda_{k+2}
},
\]
followed by an application of $Q_{k+2,k}$ to the result:
\[
\kbordermatrix{
    & k            & k+1 & k+2                      \\
k   & \ctilde_{2k} &     & \phantom{-} \stilde_{2k} \\
k+1 &              & 1   &                          \\
k+2 & \stilde_{2k} &     &          -  \ctilde_{2k}
}
\!
\kbordermatrix{
& k               &        & k+1          & k+2           \\
& \bb{\mu}_{k}    & \vrule & \phibar_{k}  &               \\
&                 & \vrule & \mubar_{k+1} &               \\
& \varepsilon_{k} & \vrule & \gamma_{k+1} & \lambda_{k+2}
}
=
\!\!\!\!\!
\kbordermatrix{
& k       &        & k+1             & k+2              \\
& \mu_{k} & \vrule & \phi_{k}        & \rho_{k}         \\
&         & \vrule & \mubar_{k+1}    &                  \\
&         & \vrule & \gammahat_{k+1} & \lambdahat_{k+2}
}.
\]
The reflections $\Qtilde_{k+1,k}$ and $\Qtilde_{k+2,k}$ zero $\gamma_k$ and $\varepsilon_{k}$ on the subdiagonals of $N_k$:
\begin{subequations}
  \begin{alignat}{4}
    \bb{\mu}_{k} &= \sqrt{\mubar_{k}^2 + \gammahat_{k}^2},
    \quad
    \ctilde_{2k-1} &=& \mubar_{k} / \bb{\mu}_{k},
    \quad
    \stilde_{2k-1} &= \gammahat_{k} / \bb{\mu}_{k},
    \quad
    k \le \ell-1,\\
    \mu_{k} &= \sqrt{\bb{\mu}_{k}^2 + \varepsilon_{k}^2},
    \quad
    \ctilde_{2k} &=& \bb{\mu}_{k} / \mu_{k},
    \quad
    \stilde_{2k} &= \varepsilon_{k} / \mu_{k},
    \quad
    k \le \ell-2,
  \end{alignat}
\end{subequations}
and they yield the recursion
\begin{subequations}
  \label{qr2_formula}
  \begin{alignat}{2}
    \phibar_k        &= \phantom{-}\stilde_{2k-1} \lambdahat_{k+1}, \quad &1 \le k \le \ell-1, \\
    \mubar_{k+1}     &=          - \ctilde_{2k-1} \lambdahat_{k+1}, \quad &1 \le k \le \ell-1,\\
    \phi_{k}         &= \ctilde_{2k} \phibar_k + \stilde_{2k} \gamma_{k+1}, \quad &1 \le k \le \ell-2, \\
    \gammahat_{k+1}  &= \stilde_{2k} \phibar_k - \ctilde_{2k} \gamma_{k+1}, \quad &1 \le k \le \ell-2, \\
    \rho_{k}         &= \phantom{-}\stilde_{2k} \lambda_{k+2}, \quad &1 \le k \le \ell-2, \\
    \lambdahat_{k+2} &=          - \ctilde_{2k} \lambda_{k+2}, \quad &1 \le k \le \ell-2, \\
    \mu_{\ell-1} &= \bb{\mu}_{\ell-1}, \\
    \phi_{\ell-1} &= \phibar_{\ell-1}, \\
    \mu_{\ell} &= \mubar_{\ell}.
  \end{alignat}
\end{subequations}
From~\eqref{eq:Ar_variants} and~\eqref{eq:qr_Nk}
we have
\begin{equation}
  \label{eq:Ar4}
  \|Ar_k\| = \Norm{N_k R_k y_k - \beta_1 \alpha_1 e_1 - \beta_1 \beta_2 e_2} = \Norm{\bmat{U_k \\ 0} R_k y_k - \bar{z}_k},
\end{equation}
where $\bar{z}_k := \Qtilde_k^T (\beta_1 \alpha_1 e_1 + \beta_1 \beta_2 e_2) = (z_k, \bb{\zeta}_{k+1}, \zetabar_{k+2})$, $k \le \ell-2$, $z_k = (\zeta_1, \dots, \zeta_{k})$ represents the first $k$ components of $\bar{z}_k$, and the recurrence starts with $\bar{z}_0 := (\bb{\zeta}_1, \zetabar_2) = (\beta_1 \alpha_1, \beta_1 \beta_2)$.
We can 
determine $\bar{z}_k$ from $\bar{z}_{k-1}$ because $\bar{z}_k = \Qtilde_{k+2,k} \Qtilde_{k+1,k} (\bar{z}_{k-1}, 0)$ for $k \le \ell-2$:
\[
\kbordermatrix{
    & k            & k+1 & k+2                      \\
k   & \ctilde_{2k} &     & \phantom{-} \stilde_{2k} \\
k+1 &              & 1   &                          \\
k+2 & \stilde_{2k} &     &          -  \ctilde_{2k}
}
\!
\kbordermatrix{
& k              & k+1                        & k+2 \\
& \ctilde_{2k-1} & \phantom{-} \stilde_{2k-1} &     \\
& \stilde_{2k-1} &          -  \ctilde_{2k-1} &     \\
&                &                            & 1
}
\!
\kbordermatrix{
&                \\
& \bb{\zeta}_k   \\
& \zetabar_{k+1} \\
& 0
}
=
\!\!\!\!\!
\kbordermatrix{
& \\
& \zeta_k          \\
& \bb{\zeta}_{k+1} \\
& \zetabar_{k+2}
}
,\]
and $z_{\ell} = z_{\ell-1} = \Qtilde_{\ell, \ell-1} \bar{z}_{\ell-2}$.
The elements are updated according to
\begin{subequations}
  \label{eq:zk}
  \begin{alignat}{2}
  \mr{\zeta}_{k}   &= \ctilde_{2k-1} \bb{\zeta}_k + \stilde_{2k-1} \zetabar_{k+1}, \quad &k \le \ell-1,\\
  \bb{\zeta}_{k+1} &= \stilde_{2k-1} \bb{\zeta}_k - \ctilde_{2k-1} \zetabar_{k+1}, \quad &k \le \ell-1,\\
  \zeta_{k}        &= \ctilde_{2k} \mr{\zeta}_{k}, \quad &k \le \ell-2, \\
  \zetabar_{k+2}   &= \stilde_{2k} \mr{\zeta}_{k}, \quad &k \le \ell-2, \\
  \zeta_{\ell-1}   &= \mathring{\zeta}_{\ell-1},&\\
  \zeta_{\ell}     &= \bb{\zeta}_{\ell}.&
  \end{alignat}
\end{subequations}
For $k \le \ell-1$, $U_k$ and $R_k$ are nonsingular, and from~\eqref{eq:Ar4}, $\|Ar_k\|$ is minimized when $U_k R_k y_k = z_k$, giving
\begin{equation}
    \label{eq:Ar_estimate}
    \norm{Ar_k} = \sqrt{\bb{\zeta}_{k+1}^2 + \bar{\zeta}_{k+2}^2},
    \quad k \le \ell-2,
    \quad \norm{Ar_{\ell-1}} = |\zeta_{\ell}|.
\end{equation}

\subsection{Computation of \texorpdfstring{$x_k$}{xk}}

Suppose $R_k$ and $U_k$ are nonsingular.
If we were to update $x_k$ directly from $x_k = V_k y_k$, all components of $y_k$ would have to be recomputed because of the backward substitutions required to solve $U_k R_k y_k = z_k$, which would require us to store $V_k$ entirely.
To avoid such drawbacks, we employ the strategy of \citet{paige-saunders-1975}.
Thus, we define $W_k$ and $D_k$ by the lower triangular systems
$R_k^T W_k^T = V_k^T$ and $U_k^T D_k^T = W_k^T$.
Then
\begin{equation}
  \label{eq:xk}
  x_k = V_k y_k = W_k R_k y_k = D_k U_k R_k y_k = D_k z_k.
\end{equation}
The columns of $W_k$ and $D_k$ are obtained from the recursions
\begin{align*}
   w_1 & = v_1 / \lambda_1, \qquad w_2 = (v_2 - \gamma_1 w_1) / \lambda_2,
\\ w_k & = (v_k - \gamma_{k-1} w_{k-1} - \varepsilon_{k-2} w_{k-2}) / \lambda_k, \quad k \ge 3,
\\ d_1 & = w_1 / \mu_1, \qquad \! d_2 = (w_2 - \phi_1 d_1) / \mu_2,
\\ d_k & = (w_k - \phi_{k-1} d_{k-1} - \rho_{k-2} d_{k-2}) / \mu_k, \quad k \ge 3,
\end{align*}
and the solution $x_k = D_k z_k$ may be updated efficiently via $x_0 = 0$ and
  \begin{equation}
    \label{eq:update_solution}
    x_k = x_{k-1} + \zeta_k d_k.
  \end{equation}
This is possible for all $k \le \ell$ if $Ax = b$ is consistent, and $k \le \ell-1$ otherwise. 
If $Ax = b$ is consistent, from~\Cref{theorem:minlength_consistent}, the final \MINARES iterate $x_{\ell}$ satisfies $r_{\ell} = 0$ and is the minimum-length solution.
If $Ax = b$ is inconsistent, from~\Cref{theorem:zetal}, $Ar_{\ell-1} = 0$.
We obtain a solution $x$ that satisfies $A^2 x = Ab$ in both cases.

\begin{theoremE}[][end, restate]
  \label{theorem:minlength_consistent}
  If $b \in \text{range}(A)$, the final \MINARES iterate $x_{\ell}$ is the minimum-length solution of $Ax = b$ (and $r_{\ell} = b - Ax_{\ell} = 0$).
\end{theoremE}

\begin{proofEE}
The final \MINARES subproblem is $T_{\ell}^2 y_{\ell} = \beta_1 \alpha_1 e_1 + \beta_1 \beta_2 e_2$  $= T_{\ell} \beta_1 e_1$. Because $b \in \text{range}(A)$, $T_{\ell}$ is nonsingular, and the latter system is equivalent to $T_{\ell} y_{\ell} = \beta_1 e_1$, the subproblem solved by \MINRES and \MINRESQLP.
The final iterate generated by these methods is the minimum-length solution of $Ax = b$ \citep[sec.~3.2 theorem 3.1]{choi-paige-minres-2011}.
\end{proofEE}

\begin{theoremE}[][end, restate]
  \label{theorem:zetal}
  If $Ax = b$ is inconsistent, $\zeta_{\ell} = 0$ and $Ar_{\ell-1} = 0$.
\end{theoremE}

\begin{proofEE}
From \eqref{eq:Nk3}, \eqref{eq:qr_Nk} and \Cref{theorem:Ar=0}:
\begin{equation*}
z_{\ell} = \Qtilde_{\ell}\T (\beta_1 \alpha_1 e_1 + \beta_1 \beta_2 e_2) = \Qtilde_{\ell}\T T_{\ell}^2 y_{\ell} = \Qtilde_{\ell}\T N_{\ell} R_{\ell} y_{\ell} = U_{\ell} R_{\ell} y_{\ell}.
\end{equation*}
When $Ax = b$ is inconsistent, $T_{\ell}$ has rank $\ell-1$ and $r_{\ell\ell} = 0$.
Because $R_{\ell}$ and $U_{\ell}$ are upper triangular matrices, $\zeta_{\ell} = u_{\ell\ell} r_{\ell\ell} \upsilon_{\ell} = 0$, where $\upsilon_{\ell}$ is the last component of $y_{\ell}$.
From~\eqref{eq:Ar_estimate}, $Ar_{\ell-1} = 0$ when $\zeta_{\ell} = 0$.
\end{proofEE}

If the minimum-norm solution is not required, such as problems~\eqref{eq:nullvector}--\eqref{eq:singularvector}, we can
stop with $x_{\ell-1}$ and avoid the computation of $x_{\ell} = \xstar$.
We can also stop with $x_{\ell-1}$ if a preconditioner is used because the minimum-norm solution is determined in a non-Euclidean norm.

We summarize the complete procedure as \Cref{alg:minares}.

  \begin{algorithm}[t]
    \caption{%
      \MINARES
    }
    \label{alg:minares}
    \begin{algorithmic}[1]
      \Require $A$, $b$, $\epsilon_{r} > 0$, $\epsilon_{Ar} > 0$, $k_{\max} > 0$
      \State $k = 0$, \,
             $x_0 = 0$
      \State $w_{-1} = w_0 = 0$, \,
             $d_{-1} = d_0 = 0$
      \State $\varepsilon_{-1} = \varepsilon_{0} = \gamma_{0} = 0$, \,
             $\rho_{-1} = \rho_{0} = \phi_{0} = 0$
      \State $\beta_1 v_1 = b$, \,
             $q_1 = Av_1$, \,
             $\alpha_1 = v_1\T q_1$ 
      \State $q_1 = q_1 - \alpha_1 v_1$, \,
             $\beta_2 v_2 = q_1$ 
      \State $\bb{\zeta}_1 = \beta_1 \alpha_1$, \,
             $\zetabar_2 = \beta_1 \beta_2$
      \State $\chibar_1 = \beta_1$, \,
             $\lambdabar_1 = \alpha_1$, \,
             $\gammabar_1 = \beta_2$
      \State $\|r_0\| = \chibar_1$, \,
             $\|Ar_0\| = (\bb{\zeta}_{1}^2 + \bb{\zeta}_{2}^2)^{\frac12}$
      \While{$\norm{r_k} > \epsilon_{r}$~\textbf{and}~
             $\norm{Ar_k} > \epsilon_{Ar}$~\textbf{and}~
             $k \le k_{\max}$}
        \State $k \leftarrow k + 1$
        \State $q_{k+1} = A v_{k+1} - \beta_{k+1} v_{k}$, \,
               $\alpha_{k+1} = v_{k+1}^T q_{k+1}$ 
        \State $q_{k+1} = q_{k+1} - \alpha_{k+1} v_{k+1}$, \,
               $\beta_{k+2} v_{k+2} = q_{k+1}$ 
        \State $\lambda_k = (\lambdabar_k^2 + \beta_{k+1}^2)^{\frac12}$, \,
               $c_k = \lambdabar_k / \lambda_k$, \,
               $s_k = \beta_{k+1} / \lambda_k$ 
        \State $\gamma_k ~~~= c_k \gammabar_k + s_k \alpha_{k+1}$, \,
               $\varepsilon_k ~~~= \phantom{-} s_k \beta_{k+2}$
        \State  $\lambdabar_{k+1} = s_k \gammabar_{k} - c_k \alpha_{k+1}$, \,
                $\gammabar_{k+1} = -c_k \beta_{k+2}$
        \If{$k == 1$}
            \State $\mubar_k = \lambda_k$, \, $\gammahat_k = \gamma_k$ 
        \Else
            \If{$k == 2$}
                \State $\lambdahat_k = \lambda_k$
            \Else
                \State $\rho_{k-2}     =  \stilde_{2k-4} \lambda_{k}$, \, 
                       $\lambdahat_{k} = -\ctilde_{2k-4} \lambda_{k}$
            \EndIf
            \State $\phibar_{k-1} = \phantom{-} \stilde_{2k-3} \lambdahat_{k}$, \,
                   $\phi_{k-1} = \ctilde_{2k-2} \phibar_{k-1} + \stilde_{2k-2} \gamma_{k}$
            \State $\mubar_{k} ~~~= -\ctilde_{2k-3} \lambdahat_{k}$, \,
                   $\gammahat_{k} \!~~~~= \stilde_{2k-2} \phibar_{k-1} - \ctilde_{2k-2} \gamma_{k}$
        \EndIf
        \State $\bb{\mu}_k = (\mubar_k^2 + \gammahat_{k}^2)^{\frac12}$\!, \,
               $\ctilde_{2k-1} = \mubar_k / \bb{\mu}_k$, \,
               $\stilde_{2k-1} = \gammahat_k / \bb{\mu}_k$
        \State $\mu_k = (\bb{\mu}_k^2 + \varepsilon_k^2)^{\frac12}$, \,
               $\ctilde_{2k} ~~~= \bb{\mu}_k / \mu_k$, \,
               $\stilde_{2k} ~~~= \varepsilon_k / \mu_k$
        \State $\mr{\zeta}_k ~~~= \ctilde_{2k-1} \bb{\zeta}_k + \stilde_{2k-1} \zetabar_{k+1}$, \,
               $\zeta_k ~~~= \ctilde_{2k} \mr{\zeta}_k$ 
        \State $\bb{\zeta}_{k+1} = \stilde_{2k-1} \bb{\zeta}_k - \ctilde_{2k-1} \zetabar_{k+1}$, \,
               $\zetabar_{k+2} = \stilde_{2k} \mr{\zeta}_k$
        \State $w_k = (v_k - \gamma_{k-1} w_{k-1} - \varepsilon_{k-2} w_{k-2}) / \lambda_k$ 
        \State $d_k = (w_k - \phi_{k-1} d_{k-1} - \rho_{k-2} d_{k-2}) / \mu_k$ 
        \State $x_k = x_{k-1} + \zeta_k d_k$ 
        \State $\norm{Ar_k} = (\bb{\zeta}_{k+1}^2 + \bar{\zeta}_{k+2}^2)^{\frac12}$
      \algstore{minares}
    \end{algorithmic}
  \end{algorithm}

  \begin{algorithm}[t]                   
    \begin{algorithmic}[1]
    \algrestore{minares}
        \State $\chi_k = c_k \chibar_{k}$, \,
               $\chibar_{k+1} = s_k \chibar_{k}$
        \If{$k == 1$}  
            \State $\psibar_k = \mu_k$, \,
                   $\bb{\pi}_{k-1} = 0$, \,
                   $\pibar_{k} = \chi_k$
            \State $\xi_k = \zeta_k$, \,
                   $\bb{\tau}_{k-1} = 0$, \,
                   $\taubar_{k} = \xi_k / \psibar_{k}$
        \ElsIf{$k == 2$}
            \State $\bb{\psi}_{k-1} = (\psibar_{k-1}^2 + \phi_{k-1}^2)^{\frac12}$, \,
                   $\chat_{k-1} = \psibar_{k-1} / \bb{\psi}_{k-1}$, \,
                   $\shat_{k-1} = \phi_{k-1} / \bb{\psi}_{k-1}$
            \State $\thetabar_{k-1} = \shat_{2k-3} \mu_k$, \,
                   $\psibar_{k} = -\chat_{2k-3} \mu_k$
            \State $\bb{\pi}_{k-1} = \chat_{2k-3} \pibar_{k-1} + \shat_{2k-3} \chi_{k}$, \,
                   $\pibar_{k} = \shat_{2k-3} \pibar_{k-1} - \chat_{2k-3} \chi_{k}$
            \State $\xi_k = \zeta_k$, \,
                   $\bb{\tau}_{k-1} = \xi_{k-1} / \bb{\psi}_{k-1}$, \,
                   $\taubar_{k} = (\xi_{k} - \thetabar_{k-1}\taubar_{k-1}) / \psibar_{k}$
        \Else
            \State $\psi_{k-2} = (\bb{\psi}_{k-2}^2 + \rho_{k-2}^2)^{\frac12}$, \,
                   $\chat_{2k-4} = \bb{\psi}_{k-2} / \psi_{k-2}$, \,
                   $\shat_{2k-4} = \rho_{k-2} / \psi_{k-2}$
            \State $\bb{\psi}_{k-1} = (\psibar_{k-1}^2 + \delta_{k}^2)^{\frac12}$, \,
                   \!~~~~$\chat_{2k-3} = \psibar_{k-1} / \bb{\psi}_{k-1}$, \,
                   $\shat_{2k-3} = \delta_{k} / \bb{\psi}_{k-1}$
            \State $\theta_{k-2} = \chat_{2k-4} \thetabar_{k-2} + \shat_{2k-4} \phi_{k-1}$, \,
                   $\omega_{k-2} = \phantom{-} \shat_{2k-4} \mu_k$
            \State $\delta_{k} ~~~= \shat_{2k-4} \thetabar_{k-2} - \chat_{2k-4} \phi_{k-1}$, \,
                   $\eta_k ~~~~= -\chat_{2k-4} \mu_k$
            \State $\thetabar_{k-1} = \shat_{2k-3} \eta_k$, \,
                   $\psibar_{k} = -\chat_{2k-3} \eta_k$, \,
                   $\upsilon_{k} = \shat_{2k-4} \bb{\pi}_{k-2} - \chat_{2k-4} \chi_k$
            \State $\bb{\pi}_{k-1} = \chat_{2k-3} \pibar_{k-1} + \shat_{2k-3} \upsilon_{k}$, \,
                   $\pibar_{k} = \shat_{2k-3} \pibar_{k-1} - \chat_{2k-3} \upsilon_{k}$
            \State $\tau_{k-2} = \bb{\tau}_{k-2} \bb{\psi}_{k-2} / \psi_{k-2}$, \,
                   $\xi_k = \zeta_{k} - \omega_{k-2} \tau_{k-2}$
            \State $\taubar_{k-1} = (\xi_{k-1} - \theta_{k-2}\tau_{k-2}) / \bb{\psi}_{k-1}$, \,
                   $\taubar_{k} = (\xi_k - \thetabar_{k-1}\taubar_{k-1}) / \psibar_{k}$
        \EndIf
        \State $\norm{r_k} = ((\bb{\pi}_{k-1} - \bb{\tau}_{k-1})^2 + (\pibar_{k} - \taubar_{k})^2 + \chibar_{k+1}^2)^{\frac12}$
      \EndWhile
    \end{algorithmic}
\end{algorithm}



\section{Stopping rules}

The end of \Cref{alg:minares} shows how $\norm{r_k}$ and $\norm{Ar_k}$ are estimated.
They are needed for use within stopping rules.
The required norm estimates are derived next.

\subsection{Estimating \texorpdfstring{$\norm{r_k}$}{|rk|}}
  
To compute $\|r_k\|$, we need an LQ factorization
\begin{equation}
  \label{eq:lq_Uk}
  U_k = \Lhat_k \Phat_k,
  \quad
  \Lhat_k = \smallbmat{
    \psi_1
 \\ \theta_1 & \psi_2
 \\ \omega_1 & \theta_2 & \psi_3
 \\          & \ddots   & \ddots & \ddots
 \\          &          & \ddots & \ddots & \psi_{k-2}
 \\          &          &        & \ddots & \theta_{k-2} & \bb{\psi}_{k-1}
 \\ \phantom{\ddots} &  &        &        & \omega_{k-2} & \thetabar_{k-1} & \psibar_{k}
  },
\end{equation}
where $\Phat_1^T = I$, $\Phat_2^T = \Phat_{1,2}$, and $\Phat_k^T = \Phat_{k-1}^T \Phat_{k-2,k} \Phat_{k-1,k}$ ($k \ge 3$) are orthogonal.
Note that $\Lhat_k$ is the L factor of a QLP decomposition of $N_k$.
If we initialize \(\psibar_1 := \mu_1\), $\Phat_{1,2}$ is defined to zero $\phi_1$:
\[
\bmat{
\psibar_{1} & \phi_{1} \\
            & \mu_{2}    
}
\bmat{
\chat_{1} & \phantom{-} \shat_{1} \\
\shat_{1} &          -  \chat_{1}
}
=
\bmat{
\bb{\psi}_{1} &             \\
\thetabar_{1} & \psibar_{2}
},
\]
where
\begin{equation}
  \bb{\psi}_{1} = \sqrt{\psibar_{1}^2 + \phi_{1}^2},
  \quad
  \chat_{1} = \psibar_{1} / \bb{\psi}_{1},
  \quad
  \shat_{1} = \phi_{1} / \bb{\psi}_{1},
  \quad
  \thetabar_{1} = \shat_{1} \mu_2,
  \quad
  \psibar_{2} = -\chat_{1} \mu_2.
\end{equation}
For $k \ge 3$, individual factorization steps may be represented as an application of $\Phat_{k-2,k}$ to $U_{k}\Phat_{k-1}^T$:
\[
\kbordermatrix{
    & k-2             & k-1           & k          \\
k-2 & \bb{\psi}_{k-2} &               & \rho_{k-2} \\
k-1 & \thetabar_{k-2} & \psibar_{k-1} & \phi_{k-1} \\
k   &                 &               & \mu_{k}
}
\!
\kbordermatrix{
& k-2          & k-1 & k                        \\
& \chat_{2k-4} &     & \phantom{-} \shat_{2k-4} \\
&              & 1   &                          \\
& \shat_{2k-4} &     &          - \chat_{2k-4}
}
=
\!\!\!\!\!
\kbordermatrix{
& k-2          & k-1           & k          \\
& \psi_{k-2}   &               &            \\
& \theta_{k-2} & \psibar_{k-1} & \delta_{k} \\
& \omega_{k-2} &               & \eta_{k}
},
\]
followed by an application of $\Phat_{k-1,k}$ to the result:
\[
\kbordermatrix{
    & k-2          & k-1           & k          \\
k-2 & \psi_{k-2}   &               &            \\
k-1 & \theta_{k-2} & \psibar_{k-1} & \delta_{k} \\
k   & \omega_{k-2} &               & \eta_{k}
}
\!
\kbordermatrix{
& k-2 & k-1          & k                        \\
& 1   &              &                          \\
&     & \chat_{2k-3} & \phantom{-} \shat_{2k-3} \\
&     & \shat_{2k-3} &          -  \chat_{2k-3}
}
=
\!\!\!\!\!
\kbordermatrix{
& k-2          & k-1             & k           \\
& \psi_{k-2}   &                 &             \\
& \theta_{k-2} & \bb{\psi}_{k-1} &             \\
& \omega_{k-2} & \thetabar_{k-1} & \psibar_{k}
}.
\]
The reflections $\Phat_{k-2,k}$ and $\Phat_{k-1,k}$ zero $\rho_{k-2}$ and $\delta_{k}$ on the superdiagonals of $U_k$:
\begin{subequations}
  \begin{alignat}{3}
    \psi_{k-2} &= \sqrt{\bb{\psi}_{k-2}^2 + \rho_{k-2}^2},
    \quad &
    \chat_{2k-4} = \bb{\psi}_{k-2} / \psi_{k-2},
    \quad &
    \shat_{2k-4} = \rho_{k-2} / \psi_{k-2},
  \\ \bb{\psi}_{k-1} &= \sqrt{\psibar_{k-1}^2 + \delta_{k}^2},
    \quad &
    \chat_{2k-3} = \psibar_{k-1} / \bb{\psi}_{k-1},
    \quad &
    \shat_{2k-3} = \delta_{k} / \bb{\psi}_{k-1},
  \end{alignat}
\end{subequations}
and for $k \ge 3$ they yield the recursion
\begin{subequations}
  \label{lq_formula}
  \begin{alignat}{2}
     \theta_{k-2}    &= \chat_{2k-4} \thetabar_{k-2} + \shat_{2k-4} \phi_{k-1}, 
  \\ \delta_{k}      &= \shat_{2k-4} \thetabar_{k-2} - \chat_{2k-4} \phi_{k-1}, 
  \\ \omega_{k-2}    &= \phantom{-} \shat_{2k-4} \mu_k,  
  \\ \eta_k          &=          -  \chat_{2k-4} \mu_k,  
  \\ \thetabar_{k-1} &= \phantom{-} \shat_{2k-3} \eta_k, 
  \\  \psibar_{k}    &=          -  \chat_{2k-3} \eta_k. 
  \end{alignat}
\end{subequations}
Assuming orthonormality of $V_{k+1}$, we have
  \begin{align}
      \norm{r_k}  = \norm{\beta_1 e_1 - T_{k+1,k} y_k}
                 &= \left\| Q_k\T \beta_1 e_1 - \bmat{R_k \\ 0} y_k \right\|                                 \nonumber \\
                 &= \left\| \bmat{\Phat_k & \\ & 1} Q_k^T \beta_1 e_1 - \bmat{\Phat_k R_k y_k \\ 0} \right\| \nonumber \\
                 &= \left\|p_{k+1} - \bmat{t_k \\ 0}\right\|, \label{eq:norm_r}
  \end{align}
where
\begin{subequations}
  \begin{align}
     &(\chi_{1}, \dots, \chi_{k}, \chibar_{k+1}) := Q_k^T \beta_1 e_1,
  \\ &p_{k+1} := (\pi_1, \dots, \pi_{k-2}, \bb{\pi}_{k-1}, \pibar_{k}, \chibar_{k+1})
                = \bmat{\Phat_k & \\ & 1} Q_k^T \beta_1 e_1,
  \\ &t_k := (\tau_1, \dots, \tau_{k-2}, \bb{\tau}_{k-1}, \bar{\tau}_k) \quad \text{solves} \quad  \Lhat_{k} t_k = z_k.
  \end{align}
\end{subequations}
The components of $Q_k^T \beta_1 e_1$ can be updated with the relations
\begin{equation}
    \label{eq:chi_k}
        \chibar_{1} = \beta_1, \qquad
        \chi_{k} = c_{k}\chibar_{k}, \qquad
        \chibar_{k+1} = s_{k}\chibar_{k},
\end{equation}
the components of $p_{k+1}$ are updated with
\begin{subequations}
    \label{eq:pi_k}
    \begin{align}
      \pibar_{1}     &= \chi_1, \\
      \upsilon_{2}   &= \chi_2, \\
      \pi_{k-2}      &= \chat_{2k-4} \bb{\pi}_{k-2} + \shat_{2k-4} \chi_k, \quad k \ge 3, \\
      \upsilon_{k}   &= \shat_{2k-4} \bb{\pi}_{k-2} - \chat_{2k-4} \chi_k, \quad k \ge 3,\\
      \bb{\pi}_{k-1} &= \chat_{2k-3} \pibar_{k-1} + \shat_{2k-3} \upsilon_{k}, \quad k \ge 2, \\
      \pibar_{k}     &= \shat_{2k-3} \pibar_{k-1} - \chat_{2k-3} \upsilon_{k}, \quad k \ge 2,
    \end{align}
\end{subequations}
and with $\omega_{-1} = \omega_{0} = \theta_{0} = \thetabar_{0} = 0$ the components of $t_k$ are updated with
\begin{subequations}
    \label{eq:tau_k}
    \begin{align}
        \xi_{k}       &= \zeta_k - \omega_{k-2} \tau_{k-2}, \\
        \taubar_{k}   &= (\xi_{k} - \thetabar_{k-1}\taubar_{k-1})  / \psibar_{k}, \\
        \bb{\tau}_{k} &= (\xi_{k} - \theta_{k-1}\tau_{k-1}) / \bb{\psi}_{k}, \\
        \tau_{k}      &= \bb{\tau}_{k} \bb{\psi}_{k} / \psi_{k}.
    \end{align}
\end{subequations}
Using \Cref{lemma:residual} we can estimate $\|r_k\|$ from the last three elements of $p_{k+1}$ and the last two of $t_k$:
\begin{subequations}
  \label{eq:norm_rk}
  \begin{align}
  \norm{r_1} &= \sqrt{(\pibar_{1}^2 - \taubar_{1}^2) + \chibar_{2}^2}, \\
  \norm{r_k} &= \sqrt{(\bb{\pi}_{k-1} - \bb{\tau}_{k-1})^2 + (\pibar_{k} - \taubar_{k})^2 + \chibar_{k+1}^2}, \quad k \ge 2.
    \end{align}
\end{subequations}

\begin{lemmaE}[][end, restate]
  \label{lemma:residual}
  In~\eqref{eq:norm_r}, $\pi_i = \tau_i$ for $i=1,\dots,k-2$.
\end{lemmaE}

\begin{proofEE}
Let $L_{k-2}$ be the leading \((k-2)\)\(\times\)\((k-2)\) submatrix of $\Lhat_k$, and $J_{m,n}$ 
be the first $m$ rows of $I_n$.  Then
\begin{align*}
  L_{k-2} J_{k-2,k+1} p_{k+1}
    &= J_{k-2,k} \Lhat_{k} J_{k,k+1} \smallbmat{\Phat_k & 0 \\ 0 & 1} Q_k^T \beta_1 e_1
\\  &= J_{k-2,k} U_k J_{k,k+1} Q_k^T \beta_1 e_1
\\  &= J_{k-2,k+2} \Qtilde_{k}^T N_k J_{k,k+1} Q_k^T \beta_1 e_1
\\  &= J_{k-2,k+2} \Qtilde_{k}^T T_{k+2,k+1} Q_k J_{k,k+1}^T J_{k,k+1} Q_k^T \beta_1 e_1
\\  &= J_{k-2,k+2} \Qtilde_{k}^T T_{k+2,k+1} Q_k (I_{k+1} - e_{k+1} e_{k+1}^T)Q_k^T \beta_1 e_1
\\  &= J_{k-2,k+2} \Qtilde_{k}^T (\beta_1 \alpha_1 e_1 + \beta_1 \beta_2 - \chibar_{k+1} T_{k+2,k+1} Q_k e_{k+1})
\\  &= J_{k-2,k+2} (\zbar_k - \chibar_{k+1} \Qtilde_{k}^T T_{k+2,k+1} Q_k e_{k+1}).
\end{align*}
We now have $T_{k+2,k+1} Q_k e_{k+1} = -(\alpha_{k+1}c_k + \beta_{k+1} c_{k-1} s_k)e_{k+1} - c_{k} \beta_{k+2}e_{k+2}$.
Further, from the structure of the reflections composing $\Qtilde_{k}^T$, the first $k-2$ elements of $\Qtilde_{k}^T T_{k+2,k+1} Q_k e_{k+1}$ are zero.
Thus,
$$L_{k-2} (\pi_1, \dots, \pi_{k-2}) = z_{k-2}.$$
Because $L_{k-2}$ is always nonsingular,
\begin{align*}
L_{k-2} \smallbmat{\pi_1 - \tau_1 \\ \vdots \\ \pi_{k-2} - \tau_{k-2}} = 0 \quad \implies \quad \smallbmat{\pi_1 \\ \vdots \\ \pi_{k-2}} = \smallbmat{\tau_1 \\ \vdots \\ \tau_{k-2}}.
\end{align*}
\end{proofEE}

\subsection{Estimating \texorpdfstring{$\norm{Ar_k}$}{|Ar|}}

From \eqref{eq:Ar_estimate} we have
\begin{equation}
    \label{eq:Ar_estimate2}
    \norm{Ar_k} = \sqrt{\bb{\zeta}_{k+1}^2 + \bar{\zeta}_{k+2}^2},
    \quad k \le \ell-2,
    \quad \norm{Ar_{\ell-1}} = |\zeta_{\ell}|.
\end{equation}

\section{\texorpdfstring{\CAR}{CAR}}
\label{sec:car}

We now introduce \CAR, a conjugate direction method in the vein of \CGM and \CR of Hestenes and Stiefel \citep{hestenes-stiefel-1952, stiefel-1955} for solving $Ax = b$ when $A$ is SPD.
By design, \CAR is equivalent to \MINARES in exact arithmetic as both methods minimize the same quantities in the same subspace, and generate the same iterates.
The name \CAR stems from the property that successive A-residuals are conjugate with respect to $A$.
The three methods generate sequences of approximate solutions $x_k$ in the Krylov subspaces $\mathcal{K}_k(A,b)$
by minimizing a quadratic function $f(x)$:
\begin{alignat*}{3}
&f_{\CGM}(x) = \tfrac{1}{2} x\T A x - b\T x, \quad
&&\nabla f_{\CGM}(x) = -r, \quad
&&\nabla^2 f_{\CGM}(x) = A, \\
&f_{\CR}(x) = \tfrac{1}{2} \norm{Ax - b}^2, \quad
&&\nabla f_{\CR}(x) = -Ar, \quad
&&\nabla^2 f_{\CR}(x) = A^2, \\
&f_{\CAR}(x) = \tfrac{1}{2} \norm{A^{2\!} x - Ab}^2, \quad
&&\nabla f_{\CAR}(x) = -A^3 r, \quad
&&\nabla^2 f_{\CAR}(x) = A^4.
\end{alignat*}
Note that all three quadratic functions satisfy $A \nabla f(x) = -\nabla^2 f(x) r$, where \(r = b - Ax\).
Because \CAR minimizes $\norm{Ar_k}$ in $\mathcal{K}_k(A,b)$, it is an alternative version of \MINARES restricted to SPD $A$.
We can derive it as a descent method with exact linesearch.
From initial vectors $x_0 = 0$ and $r_0 = p_0 = b$, we update the iterates with $x_{k+1} = x_k + \alpha_k p_k$.
From the Taylor expansion, we can determine $\alpha_k$ that minimizes $f(x_k + \alpha p_k)$:
$$f(x_k + \alpha p_k) = f(x_k) + \alpha \nabla f(x_k)\T p_k + \tfrac{1}{2} \alpha^2 p_k\T \nabla^2 f(x_k) p_k, \quad
\alpha_k = -\frac{\nabla f(x_k)\T p_k}{p_k\T \nabla^2 f(x_k) p_k}.$$
Afterwards we update the residuals with $r_{k+1} = r_{k} - \alpha_k A p_k$ and the directions with $p_{k+1} = r_{k+1} - \sum_{j=0}^{k} \gamma_{k+1,j} p_j$ such that $\Span\{p_0, \dots, p_{k+1}\}$ forms a basis of $\mathcal{K}_{k+2}(A,b)$.
We could apply a Gram–Schmidt process to orthogonalize $p_{k+1}$ against all previous directions, but a more relevant approach is to $H$-conjugate them to derive a shorter recurrence, where $H = \nabla^2 f(x)$ is constant.
$H$-conjugacy also ensures that the vectors are linearly independent.
For $i=0, \dots, k$, $p_i\T H p_{k+1} = 0$ implies $\gamma_{k+1,i} = p_i\T H r_{k+1} / p_i\T H p_i$.
Let $\mathcal{P}_k := \Span\{p_0, \dots, p_k\} = \Span\{r_0, \dots, r_k\}$.
The exact linesearch property yields \(\nabla f(x_{k+1})^T p_k\) but also $\nabla f(x_{k+1}) \perp \mathcal{P}_k$ --- see, e.g., \citep[proof of Theorem~\(5.2\)]{nocedal-wright-2006}.
Because $Ap_i = (r_i - r_{i+1}) / \alpha_i \in \Span\{r_i, r_{i+1}\} \subset \mathcal{P}_k$ for $i = 0, \dots, k-1$, we have $p_i\T A \nabla f(x_{k+1}) = - p_i\T \nabla^2 f(x_{k+1}) r_{k+1} = -p_i\T H r_{k+1} = 0$ and $\gamma_{k+1,i} = 0$.
With $\beta_k = - \gamma_{k+1,k} = - p_k\T H r_{k+1} / p_k\T H p_k$, we obtain $p_{k+1} = r_{k+1} + \beta_k p_k$.
\begin{theoremE}[][end, restate]
  \label{theorem:alpha_beta}
  For \CGM, \CR and \CAR, we have:
  $$\alpha_k = \dfrac{\rho_k}{p_k\T H p_k} \quad \text{and} \quad \beta_k = \dfrac{\rho_{k+1}}{\rho_{k}} \quad \text{with} \quad \rho_k = -\nabla f(x_k)\T r_k.$$
\end{theoremE}

\begin{proofEE}
Let $\rho_k = -\nabla f(x_k)\T r_k$.
Because $p_k = r_k + \beta_{k-1} p_{k-1}$ and $\nabla f(x_k) \perp p_{k-1}$ (exact linesearch property), $\nabla f(x_k)\T p_k = \nabla f(x_k)\T r_k$.
Therefore,
$$\alpha_k = -\dfrac{\nabla f(x_k)\T p_k}{p_k\T H p_k} = -\dfrac{\nabla f(x_k)\T r_k}{p_k\T H p_k} = \dfrac{\rho_k}{p_k\T H p_k}.$$
Because the directions $p_{i}$ are H-conjugate, $p_k\T H p_k = p_k\T H (r_k + \beta_{k-1} p_{k-1}) = p_k\T H r_k$.
With the relations $H r_{i} = -A \nabla f(x_{i})$ and $Ap_k = (r_{k} - r_{k+1}) / \alpha_k$, we have:
$$\beta_k
\!=-\dfrac{p_k\T H r_{k+1}}{p_k\T H p_k}
\!=-\dfrac{p_k\T H r_{k+1}}{p_k\T H r_k}
\!=-\dfrac{\nabla f(x_{k+1})\T (r_k - r_{k+1})}{\nabla f(x_{k})\T (r_k - r_{k+1})} 
\!=\dfrac{\nabla f(x_{k+1})\T r_{k+1}}{\nabla f(x_k)\T r_k}
\!=\dfrac{\rho_{k+1}}{\rho_{k}},
$$
where we used the fact that \(\nabla f(x_{k+1})^T r_k = -r_{k+1}^T A^i r_k = \nabla f(x_k)^T r_{k+1} = 0\), ($i = 0$ for \CGM, $i = 1$ for \CR and $i = 3$ for \CAR).
\end{proofEE}

\CGM, \CR and \CAR require $A$ to be SPD because we then have $\alpha_k > 0$ until $r_k = 0$.
The formulations of \CGM (\Cref{alg:cg}), \CR (\Cref{alg:cr}) and \CAR (\Cref{alg:car}) compare the methods and suggest efficient implementations.
The vectors $s_k = Ar_k$, $q_k = Ap_k$, $t_k = As_k = A^2 r_k$ and $u_k = Aq_k = A^2p_k$ ultimately involve just one matrix-vector product with $A$ per iteration.
Properties of \CAR are summarized in 
\Cref{theorem:properties_car}.
By virtue of its equivalence to \MINARES in exact arithmetic, \CAR allows us to establish monotonicity of relevant quantities for \MINARES (\Cref{theorem:monotonic_car}) on SPD systems.
The proofs are strongly inspired by those in \citep{fong-saunders-2012, luenberger-1970} for similar properties of \CR and \MINRES.

\begin{triplealgorithm}
  \label{alg:cg_cr_car}
  \begin{minipage}[t]{0.325\textwidth}
    \setlength{\intextsep}{0pt}
    \begin{algorithm}[H]
    \caption{\CGM}
    \label{alg:cg}
    \begin{algorithmic}
        \Require $A$, $b$, $\epsilon > 0$
        \State $k=0$, $x_0 = 0$
        \State $r_0 = b$, $p_0 = r_0$
        \State $q_0 = Ap_0$
        \State \phantom{$t_0 = A s_0$, $u_0 = t_0$}
        \State $\rho_0 = r_0\T r_0$
        \While{$\norm{r_k} > \epsilon$}
          \State $\alpha_k = \rho_k / p_k\T q_k$
          \State $x_{k+1} = x_{k} + \alpha_k p_k$
          \State $r_{k+1} = r_{k} - \alpha_k q_k$
          \State \phantom{$s_{k+1} = A r_{k+1}$}
          \State \phantom{$t_{k+1} = A s_{k+1}$}
          \State $\rho_{k+1} = r_{k+1}\T r_{k+1}$
          \State $\beta_k = \rho_{k+1} / \rho_{k}$
          \State $p_{k+1} = r_{k+1} + \beta_{k} p_k$
          \State $q_{k+1} = A p_{k+1}$
          \State \phantom{$u_{k+1} = t_{k+1} + \beta_{k} u_k$}
          \State $k \leftarrow k + 1$
        \EndWhile
    \end{algorithmic}
  \end{algorithm}
  \end{minipage}
  \hfill
  \begin{minipage}[t]{0.325\textwidth}
    \setlength{\intextsep}{0pt}
    \begin{algorithm}[H]
    \caption{\CR}
    \label{alg:cr}
    \begin{algorithmic}
        \Require $A$, $b$, $\epsilon > 0$
        \State $k=0$, $x_0 = 0$
        \State $r_0 = b$, $p_0 = r_0$
        \State $s_0 = A r_0$, $q_0 = s_0$
        \State \phantom{$t_0 = A s_0$, $u_0 = t_0$}
        \State $\rho_0 = r_0\T s_0$
        \While{$\norm{r_k} > \epsilon$}
          \State $\alpha_k = \rho_k / \norm{q_k}^2$
          \State $x_{k+1} = x_{k} + \alpha_k p_k$
          \State $r_{k+1} = r_{k} - \alpha_k q_k$
          \State $s_{k+1} = A r_{k+1}$
          \State \phantom{$t_{k+1} = A s_{k+1}$}
          \State $\rho_{k+1} = r_{k+1}\T s_{k+1}$
          \State $\beta_k = \rho_{k+1} / \rho_{k}$
          \State $p_{k+1} = r_{k+1} + \beta_{k} p_k$
          \State $q_{k+1} = s_{k+1} + \beta_{k} q_k$
          \State \phantom{$u_{k+1} = t_{k+1} + \beta_{k} u_k$}
          \State $k \leftarrow k + 1$
        \EndWhile
    \end{algorithmic}
  \end{algorithm}
  \end{minipage}
  \hfill
  \begin{minipage}[t]{0.325\textwidth}
    \setlength{\intextsep}{0pt}
    \begin{algorithm}[H]
      \caption{\CAR}
      \label{alg:car}
      \begin{algorithmic}
        \Require $A$, $b$, $\epsilon > 0$
        \State $k=0$, $x_0 = 0$
        \State $r_0 = b$, $p_0 = r_0$
        \State $s_0 = A r_0$, $q_0 = s_0$
        \State $t_0 = A s_0$, $u_0 = t_0$
        \State $\rho_0 = s_0\T t_0$
        \While{$\norm{r_k} > \epsilon$}
          \State $\alpha_k = \rho_k / \norm{u_k}^2$
          \State $x_{k+1} = x_{k} + \alpha_k p_k$
          \State $r_{k+1} = r_{k} - \alpha_k q_k$
          \State $s_{k+1} = s_{k} - \alpha_k u_k$
          \State $t_{k+1} = A s_{k+1}$
          \State $\rho_{k+1} = s_{k+1}\T t_{k+1}$
          \State $\beta_k = \rho_{k+1} / \rho_{k}$
          \State $p_{k+1} = r_{k+1} + \beta_{k} p_k$
          \State $q_{k+1} = s_{k+1} + \beta_{k} q_k$
          \State $u_{k+1} = t_{k+1} + \beta_{k} u_k$
          \State $k \leftarrow k + 1$
        \EndWhile
      \end{algorithmic}
    \end{algorithm}
  \end{minipage}
\end{triplealgorithm}

\begin{lemmaE}[][end, restate]
  \label{lem:prelim_car}
  Let \(A\) be SPD.
  The following properties hold for \CAR and \MINARES for all \(k \geq 0\):
  \begin{enumerate}[(a)]
    \item\label{car-minares-x-update} \(\zeta_{k+1} d_{k+1} = \alpha_k p_k\)
    \item\label{car-s=Ar} \(s_k = A r_k\)
    \item\label{car-q=Ap} \(q_k = A p_k\)
    \item\label{car-t=As} \(t_k = A s_k\)
    \item\label{car-u=Aq} \(u_k = A q_k\).
  \end{enumerate}
\end{lemmaE}

\begin{proofEE}
  \ref{car-minares-x-update} follows by direct comparison of \Cref{alg:minares} and \Cref{alg:car}.

  \ref{car-s=Ar}--\ref{car-u=Aq} all hold by construction at \(k = 0\).
  By induction, assume that they also hold at index \(k \geq 0\).
  Then, \(s_{k+1} = s_k - \alpha_k u_k = A r_k - \alpha_k A q_k = A r_{k+1}\), which establishes~\ref{car-s=Ar}.
  The remaining properties follow similarly.
\end{proofEE}

\begin{theoremE}[][end, restate]
  \label{theorem:properties_car}
  Let \(A\) be SPD.
  For $(i, j) \in \{0, \dots, n-1\}^2$, the following properties hold for \CAR:
  \begin{enumerate}[(a)]
    \item\label{G1:a} $p_i\T A^4 p_j = 0$ ($i \ne j$)
    \item\label{G1:b} $r_i\T A^3 p_j = 0$ ($i > j$)
    \item\label{G1:c} $r_i\T A^3 r_j = 0$ ($i \ne j$)
    \item\label{G1:d} $\alpha_i \ge 0$
    \item\label{G1:e} $\beta_i \ge 0$
    \item\label{G1:f} $q_i\T u_j = p_i\T A^3 p_j \ge 0$
    \item\label{G1:g} $q_i\T q_j = p_i\T A^2 p_j \ge 0$
    \item\label{G1:h} $q_i\T p_j = p_i\T A p_j \ge 0$
    \item\label{G1:i} $p_i\T p_j \ge 0$
    \item\label{G1:j} $x_i\T p_j \ge 0$
    \item\label{G1:k} $r_i\T q_j = r_i\T A p_j \ge 0$.
  \end{enumerate}
\end{theoremE}

\begin{proofEE}
  Because $\nabla^2 f_{\CAR}(x) = A^4$, we $A^4$-conjugate the vectors $p_i$ by construction and~\ref{G1:a} is satisfied.

  Because $\nabla f_{\CAR}(x_{i}) = - A^3 r_{i}$, the exact linesearch property yields~\ref{G1:b} as in \citep[proof of Theorem~\(5.2\)]{nocedal-wright-2006}.

  If $i > j$, $r_i\T A^3 r_j = r_i\T A^3 (p_j - \beta_{j-1} p_{j-1}) = 0$ by~\ref{G1:b}.
  If \(i < j\), $r_i\T A^3 r_j = (p_i - \beta_{i-1} p_{i-1})\T A^3 r_j = 0$, again thanks to~\ref{G1:b}, which proves~\ref{G1:c}.

  First note that $\rho_i = s_i\T t_i = r_i\T A^3 r_i \ge 0$ because $A$ is SPD.
  Thus $\alpha_i = \rho_i / \norm{u_i}^2 \ge 0$ and $\beta_i = \rho_{i+1} / \rho_{i} \ge 0$, which proves~\ref{G1:d} and~\ref{G1:e}.

  We now establish~\ref{G1:f} by induction.
  If $i = j$, $q_i\T u_i = q_i\T A q_i \ge 0$ because $A$ is SPD.
  Assuming $q_i\T u_j \ge 0$ when $|i-j| = k-1 \ge 0$, we want to show the result for $|i-j| = k$.
  If $i-j = k > 0$ then $q_i\T u_j = q_i\T u_{i-k}$.
  Otherwise we have $j-i = k > 0$ and $q_i\T u_j = q_i\T u_{i+k}$.
  \Cref{lem:prelim_car} yields

\begin{alignat*}{4}
    & q_i\T u_{i-k} &&= (s_i + \beta_{i-1} q_{i-1})\T u_{i-k}                         & \quad \quad q_i\T u_{i+k} &= q_i\T A q_{i+k}
    \\ &                        &&= s_i\T u_{i-k} + \beta_{i-1} q_{i-1}\T u_{i-k}     & \quad \quad               &= q_i\T A (s_{i+k} + \beta_{i+k-1} q_{i+k-1})
    \\ &                        &&= r_i\T A^3 p_{i-k} + \beta_{i-1} q_{i-1}\T u_{i-k} & \quad \quad               &= p_i\T A^3 r_{i+k} + \beta_{i+k-1} u_i\T q_{i+k-1}
    \\ &                        &&= \beta_{i-1} q_{i-1}\T u_{i-k}                     & \quad \quad               &= \beta_{i+k-1} q_{i+k-1}\T u_i
\end{alignat*}
$\beta_{i-1} \ge 0$ and $\beta_{i+k-1} \ge 0$ by~\ref{G1:e}.
$q_{i-1}\T u_{i-k} \ge 0$ and $q_{i+k-1}\T u_i \ge 0$ by induction assumption.
Thus, $q_i\T u_j \ge 0$ for $|i-j| = k$, which completes the proof of~\ref{G1:f}.



  At termination, define $\mathcal{P} = \Span\{p_0, \dots, p_{\ell-1}\}$, $\mathcal{Q} = \Span\{q_0, \dots, q_{\ell-1}\} = A \mathcal{P}$
  and $\mathcal{U} = \Span\{u_0, \dots, u_{\ell-1}\} = A \mathcal{Q}$.
  By construction, $\mathcal{P} = \Span\{b, \dots, A^{\ell-1}b\}$, $\mathcal{Q} = \Span\{Ab, \dots, A^{\ell}b\}$ and $\mathcal{U} = \Span\{A^2 b, \dots, A^{\ell+1}b\}$.
  Again by construction, $x_{\ell} \in \mathcal{P}$, and since $r_{\ell} = 0$, we have $A x_{\ell} = b \in \mathcal{Q}$ and $A^2 x_{\ell} = Ab \in \mathcal{U}$.
  We see that $\mathcal{P} \subset \mathcal{Q} \subset \mathcal{U}$.

  \ref{G1:a} and \Cref{lem:prelim_car}~\ref{car-q=Ap}--\ref{car-u=Aq} imply that \(u_i\T u_j = 0\) for \(i \neq j\), and therefore,
  $\left\{u_k / \norm{u_k}\right\}_{k = 0, \dots, \ell-1}$ forms an orthonormal basis for $\mathcal{U}$.
  Thus, if we project $p_i$ and $q_i$ into $\mathcal{U}$, we have
  $$p_i = \sum_{k=0}^{\ell-1} \dfrac{p_i\T u_k}{u_k\T u_k} u_k \quad \text{and} \quad q_i = \sum_{k=0}^{\ell-1} \dfrac{q_i\T u_k}{u_k\T u_k} u_k.$$
  Scalar products between these vectors can be expressed as
  $$ q_i\T q_j = \sum_{k=0}^{\ell-1} \dfrac{(q_i\T u_k)(q_j\T u_k)}{\norm{u_k}^2},\! \quad p_i\T q_j = \sum_{k=0}^{\ell-1} \dfrac{(p_i\T u_k)(q_j\T u_k)}{\norm{u_k}^2} \quad \!\text{and}\! \quad p_i\T p_j = \sum_{k=0}^{\ell-1} \dfrac{(p_i\T u_k)(p_j\T u_k)}{\norm{u_k}^2}.$$
  Thus $q_i\T q_j \ge 0$ by~\ref{G1:f}, proving \ref{G1:g}.
  Because $p_i\T u_k = p_i\T A q_k = q_i\T q_k$, $p_i\T q_j \ge 0$ and $p_i\T p_j \ge 0$ by~\ref{G1:f} and~\ref{G1:g}, which proves~\ref{G1:h} and~\ref{G1:i}.

  By construction, $x_i = \sum_{k=0}^{i} \alpha_k p_k$ and so $x_i\T p_j \ge 0$ by~\ref{G1:d} and~\ref{G1:i}, proving~\ref{G1:j}.

  Finally, $r_i\T q_j = \sum_{k=i}^{\ell-1} \alpha_k q_k\T q_j \ge 0$ by~\ref{G1:d} and~\ref{G1:g}, proving~\ref{G1:k}.
\end{proofEE}

\begin{theoremE}[][end, restate]
  \label{theorem:monotonic_car}
  For \CAR (and hence \MINARES) applied to $Ax = b$ when $A$ is SPD, the following properties are satisfied:
  \begin{itemize}
    \item $\|x_k\|$ increases monotonically
    \item $\|\xstar - x_k\|$ decreases monotonically
    \item $\|\xstar - x_k\|_A$ decreases monotonically
    \item $\|r_k\|$ decreases monotonically.
  \end{itemize}
\end{theoremE}


\begin{proofEE}
  From \Cref{theorem:properties_car}~\ref{G1:d} and~\ref{G1:j},
  \begin{align*}
      \norm{x_k}^2 - \norm{x_{k-1}}^2 &= (x_{k-1} + \alpha_k p_k)\T (x_{k-1} + \alpha_k p_k) - x_{k-1}\T x_{k-1}
      \\                              &= 2 \alpha_k p_k\T x_{k-1} + \alpha_k^2 \norm{p_k}^2 \ge 0.
  \end{align*}
  From \Cref{theorem:properties_car}~\ref{G1:d} and~\ref{G1:i},
  \begin{align*}
      \norm{\xstar - x_{k-1}}^2 - \norm{\xstar - x_{k}}^2 &= \left(\sum_{i=k}^{\ell-1} \alpha_{i} p_i\right)\TT \left(\sum_{i=k}^{\ell-1} \alpha_{i} p_i\right) - \left(\sum_{i=k+1}^{\ell-1} \alpha_{i} p_i\right)\TT \left(\sum_{i=k+1}^{\ell-1} \alpha_{i} p_i\right)
      \\ &= 2 \alpha_k p_k\T \left(\sum_{i=k+1}^{\ell-1} \alpha_i p_i\right) + \alpha_k^2 \norm{p_k}^2 \ge 0.
  \end{align*}
  From \Cref{theorem:properties_car}~\ref{G1:d} and~\ref{G1:h},
  \begin{align*}
      \norm{\xstar\!\! - \! x_{k-1}}_A^2 - \norm{\xstar\!\! - \! x_{k}}_A^2 \! &= \! \left(\sum_{i=k}^{\ell-1} \alpha_{i} p_i\right)\TT A \left(\sum_{i=k}^{\ell-1} \alpha_{i} p_i\right) \!-\! \left(\sum_{i=k+1}^{\ell-1} \alpha_{i} p_i\right)\TT A \left(\sum_{i=k+1}^{\ell-1} \alpha_{i} p_i\right)\!\!
  \\ &= 2 \alpha_k q_k\T \left(\sum_{i=k+1}^{\ell-1} 
  \alpha_i p_i\right) + \alpha_k^2 q_k\T p_k \ge 0.
  \end{align*}
  From \Cref{theorem:properties_car}~\ref{G1:d} and~\ref{G1:k},
  \begin{align*}
        \norm{r_{k-1}}^2 - \norm{r_k}^2 &= r_{k-1}\T r_{k-1} - r_{k}\T r_k
      \\                                &= (r_{k} + \alpha_{k-1} q_{k-1})\T (r_{k} + \alpha_{k-1} q_{k-1}) - r_{k}\T r_k
      \\                                &= 2 \alpha_{k-1} q_{k-1}^T r_k + \alpha_{k-1}^2 \|q_{k-1}\|^2 \ge 0.
  \end{align*}
\end{proofEE}

\section{Implementation and numerical experiments}

We implemented \Cref{alg:minares} and \Cref{alg:car}
in Julia \citep{bezanson-edelman-karpinski-shah-2017}, version \(1.9\), as part of our \texttt{Krylov.jl} collection of Krylov methods \citep{montoison-orban-2023}.
These implementations of \MINARES and \CAR are applicable in any floating-point system supported by Julia, including complex numbers, and they run on CPU and GPU.
They also support preconditioners.

We evaluate the performance of \MINARES on systems generated from symmetric matrices $A$ in the SuiteSparse Matrix Collection \citep{davis-hu-2011}.
In each case we first scale $A$ to be $A/\alpha$ with $\alpha = \max |A_{ij}|$, so that $\norm{A} \approx 1$.

\begin{figure}[t] 
  \centering
  \includetikzgraphics[width=0.49\textwidth]{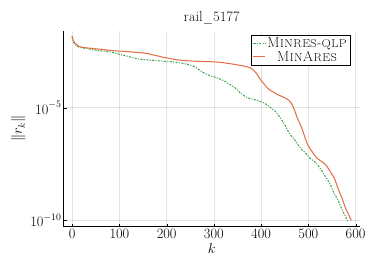}
  \hfill
  \includetikzgraphics[width=0.49\textwidth]{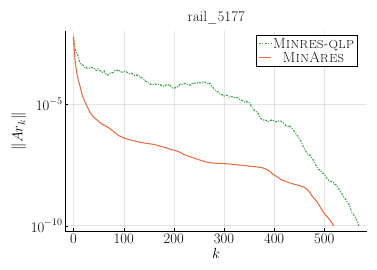}
\\
  \includetikzgraphics[width=0.49\textwidth]{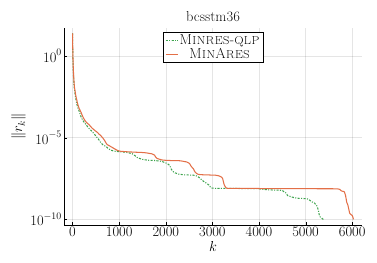}
  \hfill
  \includetikzgraphics[width=0.49\textwidth]{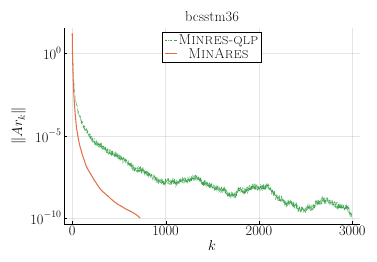}
  \vspace{-2\baselineskip}
  \caption{\label{fig:ufl1}Residual and A-residual histories for \MINARES and \MINRESQLP on consistent systems generated from the SuiteSparse Matrix Collection.
  Top: System based on the nonsingular matrix rail\_5177
  $(n = 5177)$.
  Bottom: System based on the singular matrix bcsstm36
  $(n = 23052)$.}
\end{figure}

In our first set of experiments, we compare \MINARES to our Julia implementation of \MINRESQLP in terms of number of iterations on consistent systems when the stopping criterion is $\|r_k\| \leq 10^{-10}$, then when it is $\|Ar_k\| \leq 10^{-10}$.
The right-hand side $b = Ae$ (with $e$ a vector of ones) ensures that the system is consistent even if $A$ is singular. 
The residual and A-residual are calculated explicitly at each iteration in order to evaluate $\|r_k\|$ and $\|Ar_k\|$. 
(To get a fair comparison, \eqref{eq:norm_rk} and~\eqref{eq:Ar_estimate2} are not used.)
\Cref{fig:ufl1} reports residual and A-residual histories for \MINARES and \MINRESQLP on problems \emph{rail\_5177} and \emph{bcsstm36}.
We observe that \MINRESQLP's $\norm{A r_k}$ is erratic, whereas \MINARES's $\norm{A r_k}$ and $\norm{r_k}$ are both smooth.
Also, \MINRESQLP's $\norm{Ar_k}$ lags further behind \MINARES's than \MINARES's $\norm{r_k}$ does behind \MINRESQLP's.
When the system is consistent, we have similar behavior whether $A$ is singular or not.

\begin{figure}[t]  
  \centering
  \includetikzgraphics[width=0.49\textwidth]{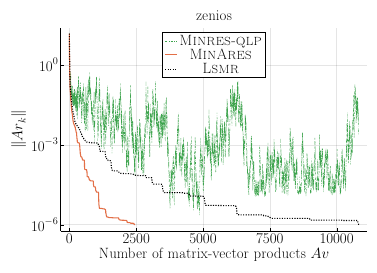}
  \hfill
  \includetikzgraphics[width=0.49\textwidth]{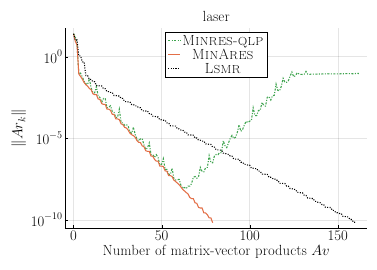}
  \vspace{-2\baselineskip}
  \caption{\label{fig:ufl2}A-residual history for \MINARES, \MINRESQLP and \LSMR on singular inconsistent systems generated from the SuiteSparse Matrix Collection.
  Left: System based on the singular matrix zenios
  $(n = 2873)$.
  Right: System based on the singular matrix laser
  $(n = 3002)$.}
\end{figure}

In a second set of experiments, we compare \MINARES to our Julia implementations of \MINRESQLP and \LSMR in terms of number of matrix-vector products $Av$ on singular inconsistent systems with $b=e$ when the stopping criterion is $\|Ar_k\| \leq 10^{-6}$ for the problem \emph{zenios} and $\|Ar_k\| \leq 10^{-10}$ for \emph{laser}.
\Cref{fig:ufl2} shows that \MINRESQLP has difficulty reaching the specified $\norm{A r_k}$,
but \MINARES performs well and converges much faster than \LSMR, the only other Krylov method that minimizes $\norm{A r_k}$.



\section{Summary}

\MINARES completes the family of Krylov methods based on the symmetric Lanczos process.
By minimizing $\norm{A r_k}$ (which always converges to zero), \MINARES can be applied safely to any symmetric system.
For SPD systems, \CAR is equivalent to \MINARES and extends the conjugate directions family \CGM and \CR.
For such systems we prove that $\norm{r_k}$, $\norm{x_k - \xstar}$ and $\norm{x_k - \xstar}_A$ decrease monotonically for \CAR and hence \MINARES.

On consistent symmetric systems, \MINARES is a relevant alternative to \MINRES and \MINRESQLP because it converges in a similar number of iterations if the stopping condition is based on $\norm{r_k}$, and much faster if the stopping condition is based on $\norm{A r_k}$.
On singular inconsistent symmetric systems, \MINARES outperforms \MINRESQLP and \LSMR, and should be the preferred method.

\subsection*{Acknowledgements}

This work began while the first author was visiting ICME at Stanford University in Spring 2022.
My thanks to Mike Saunders for making my stay such a pleasant one.



\appendix

\section{Proofs}

\printProofs

\small
\bibliographystyle{abbrvnat}
\bibliography{abbrv,minares}
\normalsize

\end{document}